\documentclass[12pt]{amsart}
\usepackage{latexsym,amsmath,amssymb,amsthm,amsfonts}
\usepackage[numeric]{amsrefs}
\usepackage{tikz}
\usetikzlibrary{arrows.meta,calc}
\usepackage{needspace}
\usepackage[section]{placeins}
\usepackage[hidelinks]{hyperref}
\usepackage[capitalise]{cleveref}
\newcommand{\R}{\mathbb R}
\newcommand{\RP}{\mathbb{RP}}
\newcommand{\Sph}{\mathbb S}
\DeclareMathOperator{\tr}{tr}
\DeclareMathOperator{\diam}{diam}

\DeclareMathOperator{\sconv}{sconv}
\DeclareMathOperator{\area}{area}
\DeclareMathOperator{\Sym}{Sym}
\newcommand{\norm}[1]{\lVert#1\rVert}
\newtheorem{theorem}{Theorem}
\newtheorem{lemma}[theorem]{Lemma}
\newtheorem{proposition}[theorem]{Proposition}
\newtheorem{corollary}[theorem]{Corollary}
\theoremstyle{remark}

\title[A diameter-one set in $\mathbb{R}^5$]{A diameter-one set in $\mathbb{R}^5$ that cannot be covered by six diameter-one balls}
\author{A. Arman}
 \address{Department of Mathematics, University of Manitoba, Winnipeg, MB, R3T 2N2, Canada}
 \email{andrew0arman@gmail.com}
 
\author{A.\ Bondarenko}
 \address{Department of Mathematical Sciences, Norwegian University of Science and
 	Technology, NO-7491 Trondheim, Norway}
 \email{andriybond@gmail.com}
 \thanks{The second author was supported in part by Grant 334466 of the Research Council of Norway.}

\author{A.\ Prymak}
 \address{Department of Mathematics, University of Manitoba, Winnipeg, MB, R3T 2N2, Canada}
 \email{prymak@gmail.com}
 \thanks{The third author was supported by NSERC of Canada Discovery Grant RGPIN-2026-06488.}

\author{D.\ Radchenko}
 \address{Institut des Hautes \'Etudes Scientifiques, CNRS, Laboratoire Alexander Grothendieck, 35 route de Chartres, Bures-sur-Yvette 91440, France}
 \email{danradchenko@gmail.com}
\thanks{The fourth author was supported by ERC Starting Grant No. 101078782.}

\date{November 5, 2026}
\hypersetup{pdftitle={A five-dimensional example for Grunbaum's covering problem},pdfauthor={A. Arman, A. Bondarenko, A. Prymak, D. Radchenko}}
\begin{document}
\begin{abstract}
We construct a compact set of diameter $1$ in $\R^5$ which cannot be covered by six Euclidean balls of diameter $1$. Its convex hull is a full-dimensional convex body satisfying the same property. Moreover, a finite subset of diameter $1$ cannot be covered by six closed balls even when their radius is increased to $0.50001$. The compact set is the quadratic image of the unit sphere in the space of symmetric $3\times3$ matrices of trace zero. A hypothetical cover leads to twelve spherical ellipses covering $\Sph^2$ and arranged in six antipodal pairs. Using Euler's formula, we replace them by twelve pentagonal regions contained in the ellipses. An area estimate then gives a contradiction.
\end{abstract}
\maketitle

\section{Introduction}
For a bounded set $A\subset\R^n$, write $\diam A$ for its Euclidean diameter and let $N_{\rm cov}(A,r)$ be the least number of closed balls of radius $r$ needed to cover $A$. Gr\"unbaum's equal-diameter covering problem~\cite{Grunbaum} asks about
\[
 g_n=\sup_{\diam A=1}N_{\rm cov}(A,1/2).
\]
For large $n$, Bourgain and Lindenstrauss~\cite{BL} proved an exponential upper bound for $g_n$, while recent exponential lower bounds appear in~\cites{ABP,ABPR}.

For comparison, let $b_n$ be the least integer such that every bounded set of positive diameter in $\R^n$ can be partitioned into at most $b_n$ subsets of strictly smaller diameter. Borsuk~\cite{Borsuk} asked whether $b_n=n+1$ for every $n$. Kahn and Kalai~\cite{KK} disproved this in sufficiently high dimensions. More recently, OpenAI~\cite{OpenAI} gave a compact example in dimension nine which cannot be covered by ten sets of strictly smaller diameter. Their example consists of the rank-one orthogonal projectors onto lines in $\R^4$, viewed in the nine-dimensional affine space of symmetric $4\times4$ matrices of trace one with the Frobenius metric.

In Gr\"unbaum's problem the covering sets must be balls, although their diameter is allowed to equal that of the original set. The restriction on their shape suggests that the analogous equality $g_n=n+1$ might fail in a lower dimension. We show that it fails already in dimension five. Our example uses the same rank-one projector construction as the nine-dimensional example, with $\R^3$ in place of $\R^4$. The projector set in~\cite{OpenAI} is infinite and is not discrete. Our argument also yields a finite counterexample to the $n+1$-ball covering assertion in dimension five.

Let $\Sym_0(3)$ be the five-dimensional Euclidean space of real symmetric $3\times3$ matrices of trace zero, with scalar product $\langle A,B\rangle=\tr(AB)$ and norm $\norm{A}^2=\tr(A^2)$. Write $I$ for the $3\times3$ identity matrix and $\Sph^2$ for the unit sphere in $\R^3$. Consider the quadratic image
\begin{equation}\label{eq:veronese}
 p(u)=\frac{uu^T-I/3}{\sqrt2},\qquad
 V=\{p(u):u\in\Sph^2\}\subset\Sym_0(3).
\end{equation}
Since $p(u)=p(-u)$, this is the classical Veronese model of the real projective plane $\RP^{2}$, see also~\cites{CHS,Kalai}.

\begin{theorem}\label{thm:main}
The compact set $V\subset\Sym_0(3)\cong\R^5$ has diameter $1$ and cannot be covered by six closed balls of radius $1/2$. Its convex hull is a full-dimensional convex body of diameter $1$ which also cannot be covered by six such balls. In particular, $g_5\ge7$.
\end{theorem}

\begin{corollary}\label{cor:finite}
There is a finite set $F\subset V$ of diameter $1$ which cannot be covered by six closed balls of radius $0.50001$.
\end{corollary}

The proof reduces the covering problem to twelve spherical ellipses. Euler's formula determines their intersections and gives twelve pentagonal regions covering $\Sph^2$. Each region lies in the spherical convex hull of five points of one ellipse, and each such hull has area less than $\pi/3$. Section~\ref{sec:topological} gives an alternative topological proof of the intersection structure in Proposition~\ref{prop:structure}. In Section~\ref{sec:Borsuk}, we partition $V$ into six sets of smaller diameter. Thus our construction gives $g_5>6$ but does not establish $b_5>6$.

\section{The quadratic image and its ball sections}\label{sec:sections}

\subsection{Distances and the convex hull}
A direct calculation gives
\begin{equation}\label{eq:distance}
 \norm{p(u)}^2=\frac13,\qquad
 \norm{p(u)-p(v)}^2=1-(u\cdot v)^2,
 \qquad u,v\in\Sph^2.
\end{equation}
These identities follow by expanding the squared Frobenius norms and using $\tr(uu^Tvv^T)=(u\cdot v)^2$. In particular, $V$ is compact by continuity and has diameter $1$, attained whenever $u\perp v$.

In orthonormal coordinates,
\begin{equation}\label{eq:coordinates}
 V=\left\{\left(\frac{x^2-y^2}{2},
 \frac{x^2+y^2-2z^2}{\sqrt{12}},xy,xz,yz\right):
 x^2+y^2+z^2=1\right\}\subset\R^5.
\end{equation}
Its convex hull is
\begin{equation}\label{eq:convex-body}
 K=\left\{\frac{Q-I/3}{\sqrt2}:Q\succeq0,\ \tr Q=1\right\}.
\end{equation}
Here $Q\succeq0$ means that $Q$ is symmetric and positive semidefinite. Formula~\eqref{eq:convex-body} follows from the spectral decomposition of $Q$ into orthogonal rank-one projections with nonnegative coefficients summing to one. The positive definite matrix $I/3$ has a neighbourhood of positive definite matrices in the trace-one affine space. Hence $K$ is a full-dimensional convex body in $\Sym_0(3)$.

Taking a convex hull preserves diameter. More explicitly, if $x=\sum_r\alpha_rx_r$ and $y=\sum_s\beta_sy_s$, where $x_r,y_s\in V$ and both collections of coefficients are nonnegative and sum to one, then
\[
 \norm{x-y}\le\sum_{r,s}\alpha_r\beta_s\norm{x_r-y_s}\le1.
\]
Thus $\diam K=1$. Since $V\subset K$, it remains to show that six balls of radius $1/2$ cannot cover $V$.

\subsection{Spherical notation and gnomonic coordinates}
The distance on $\Sph^2$ will always be angular distance,
\[
 d_{\Sph}(u,v)=\arccos(u\cdot v).
\]
For nonantipodal points $u,v$, their \emph{short arc} is the shorter great-circle arc joining them. An open hemisphere is a set $H_w=\{u\in\Sph^2:u\cdot w>0\}$, where $w\in\Sph^2$. The gnomonic map from the plane $w^\perp=\{z:z\cdot w=0\}$ onto $H_w$ is
\[
 G_w(z)=\frac{w+z}{\sqrt{1+\norm z^2}}.
\]
Geometrically, it projects the affine plane $w+w^\perp$ radially onto the sphere. Its inverse sends short arcs to line segments, because the rays through a great circle lie in a plane through the origin.

A set contained in an open hemisphere is \emph{spherically convex} if it contains the short arc between each two of its points. In gnomonic coordinates this is ordinary convexity. For a finite set in an open hemisphere, its spherical convex hull, denoted by $\sconv$, is the gnomonic image of its planar convex hull; equivalently, it consists of the normalized nonzero convex combinations of its points. A \emph{spherical ellipse} in this paper means the filled image under $G_w$ of a planar ellipse centred at zero, with positive semiaxes. All ellipses below include their boundaries.

\subsection{The inverse image of a ball}
We now describe the inverse image under $p$ of a ball of radius $1/2$.

\begin{lemma}\label{lem:sections}
Let $A\in\Sym_0(3)$. The inverse image under $p$ of the closed ball with centre $A$ and radius $1/2$ is contained in a set
\[
 E_M=\{u\in\Sph^2:u^TMu\ge\tau\},\qquad
 \tau=\frac1{\sqrt6},
\]
where $M\in\Sym_0(3)$ and $\norm M=1$. Every such $E_M$ is either a great circle or the union of two antipodal spherical ellipses. In the latter case their gnomonic semiaxes satisfy
\begin{equation}\label{eq:axes}
 a^2+h^2=1,
\end{equation}
and each ellipse has spherical diameter less than $\pi/2$.
\end{lemma}

\begin{proof}
A ball centred at zero misses $V$, by~\eqref{eq:distance}; its empty inverse image is contained in any $E_M$. Otherwise write $A=\rho M$, where $\rho=\norm A>0$, $\tr M=0$ and $\norm M=1$. Expanding the square gives
\begin{equation}\label{eq:ball-section}
 \norm{p(u)-A}^2\le\frac14
 \quad\Longleftrightarrow\quad
 u^TMu\ge\frac{\rho^2+1/12}{\sqrt2\rho}
 \ge\frac1{\sqrt6}.
\end{equation}
The last inequality follows from $\rho+(12\rho)^{-1}\ge2/\sqrt{12}$.

Let $\lambda_1\ge\lambda_2\ge\lambda_3$ be the eigenvalues of $M$. Their sum is zero and the sum of their squares is one. We first verify
\begin{equation}\label{eq:eigenvalues}
 \lambda_1\ge\tau\ge\lambda_2,
\end{equation}
with equality on either side only for $(\lambda_1,\lambda_2,\lambda_3)=(\tau,\tau,-2\tau)$. If $\lambda_2\ge0$, then $\lambda_1\ge\lambda_2$ and $\lambda_3=-\lambda_1-\lambda_2\le-2\lambda_2$, whence $1\ge6\lambda_2^2$. If $\lambda_2<0$, the upper bound on $\lambda_2$ is immediate. For the other bound put $s=\lambda_1>0$. The ordering and trace condition imply $-s/2\le\lambda_2\le s$, and
\[
 1=s^2+\lambda_2^2+(s+\lambda_2)^2\le6s^2.
\]
The equality conditions follow from the same calculations.

For the exceptional eigenvalue triple, $u^TMu\ge\tau$ holds exactly on the great circle perpendicular to a $\lambda_3$ eigenvector. In every other case $\lambda_1>\tau>\lambda_2$. Choose an orthonormal eigenbasis $e_1,e_2,e_3$ and write
\begin{equation}\label{eq:gnomonic}
 u=\frac{e_1+\xi e_2+\eta e_3}{\sqrt{1+\xi^2+\eta^2}}.
\end{equation}
There are no points of $E_M$ on $u\cdot e_1=0$. In the hemisphere $u\cdot e_1>0$, its defining inequality becomes
\[
 \beta\xi^2+\gamma\eta^2\le r,
 \qquad r=\lambda_1-\tau,\quad
 \beta=\tau-\lambda_2,\quad\gamma=\tau-\lambda_3.
\]
All three coefficients are positive. The trace and norm conditions, using $6\tau^2=1$, give
\[
 \beta+\gamma=3\tau+r,\qquad
 \beta\gamma=r(\beta+\gamma).
\]
Consequently $a^2=r/\beta$, $h^2=r/\gamma$ and $a^2+h^2=1$. The other component is its antipodal image. Both semiaxes are less than one, so every point in the upper ellipse is at angular distance at most $\arctan\max(a,h)<\pi/4$ from $e_1$. The triangle inequality gives spherical diameter less than $\pi/2$.
\end{proof}

\begin{figure}[ht]
\centering
\begin{tikzpicture}[>=Latex,font=\small]
 \begin{scope}[xshift=-3.4cm,scale=.93]
  \fill[gray!8] (0,0)--(1.65,2.15)--(3.05,2.15)--cycle;
  \draw (0,0) node[below]{$O$}--(3.05,2.15);
  \draw (0,0)--(1.65,2.15);
  \draw[gray] (-.35,1.55)--(1,2.65)--(3.85,2.65)--(2.5,1.55)--cycle;
  \draw[thick] (1.65,2.15)--(3.05,2.15);
  \draw[thick] (35.18:1.25) arc[start angle=35.18,end angle=52.5,radius=1.25];
  \node[anchor=west] at (1.25,.52) {short arc};
  \node at (2.4,2.92) {line segment};
  \node at (1.7,-.48) {radial projection};
 \end{scope}
 \begin{scope}[xshift=3.2cm,yshift=1.15cm]
  \draw[->,gray] (-1.75,0)--(1.9,0) node[right]{$\xi$};
  \draw[->,gray] (0,-1.35)--(0,1.5) node[above]{$\eta$};
  \draw[thick,fill=gray!10] (0,0) ellipse (1.5 and .94);
  \draw[dashed] (-1.08,-.65)--(.25,.92);
  \fill (-1.08,-.65) circle(1.6pt);\fill (.25,.92) circle(1.6pt);
  \node at (0,-1.62) {$\xi^2/a^2+\eta^2/h^2\le1$};
 \end{scope}
\end{tikzpicture}
\caption{Gnomonic coordinates turn short spherical arcs into line segments. A component of $E_M$ becomes a planar ellipse with $a^2+h^2=1$. The left panel is a schematic view of the rays and the affine plane.}
\label{fig:gnomonic}
\end{figure}
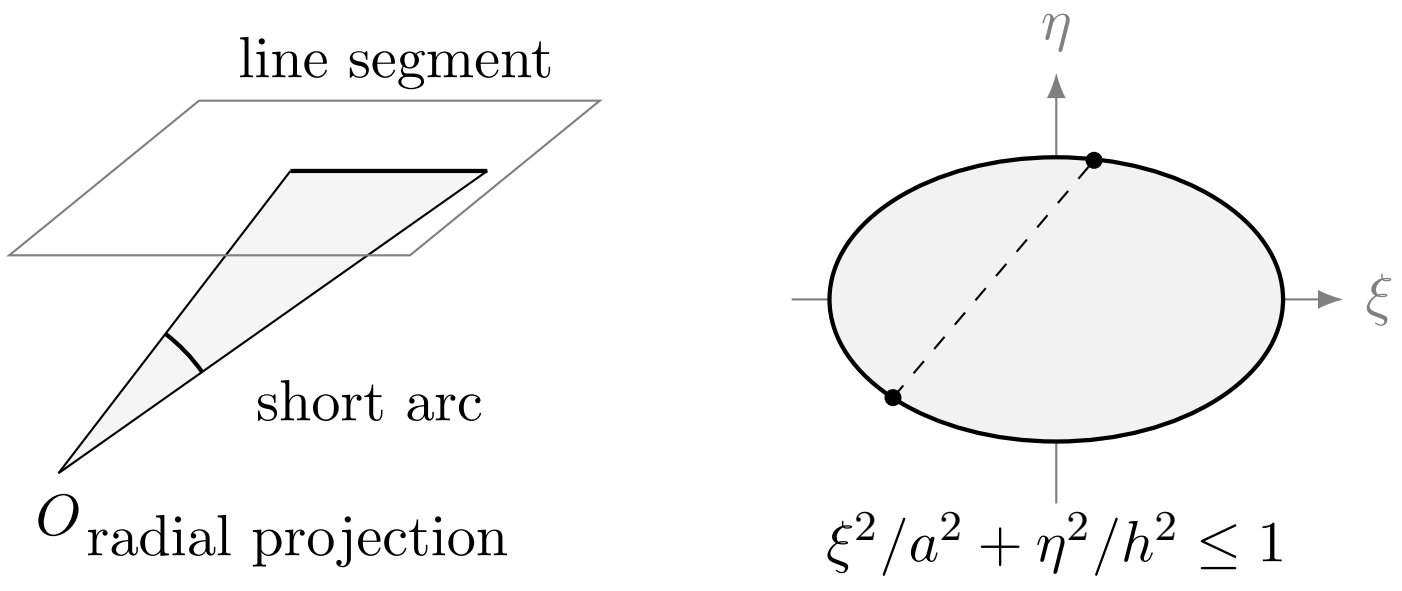

\section{Area estimates}\label{sec:area}
We use $\area$ for the usual surface area on the unit sphere, so $\area(\Sph^2)=4\pi$. The area element in~\eqref{eq:gnomonic} is
\begin{equation}\label{eq:jacobian}
 (1+\xi^2+\eta^2)^{-3/2}\,d\xi\,d\eta,
\end{equation}
which follows by differentiating the parametrization. Since the density decreases with the distance from the origin, a centred disk has the largest spherical image among planar sets of a given area.

\begin{lemma}\label{lem:area}
Let $C$ be a spherical ellipse with gnomonic semiaxes $a,h$ satisfying $a^2+h^2\le1$. Then
\begin{equation}\label{eq:ellipse-bound}
 \area(C)\le B_0:=2\pi\left(1-\sqrt{\frac23}\right)<\frac{2\pi}{5}.
\end{equation}
The spherical convex hull of any five points of $C$ has area strictly less than $\pi/3$.
\end{lemma}

\begin{proof}
Work in a gnomonic chart $G_w$ for $C$. If $D$ is a planar measurable set of area $A$, let $D_0$ be the disk centred at zero with the same area. The sets $D\setminus D_0$ and $D_0\setminus D$ have equal area. On the former the density in~\eqref{eq:jacobian} is at most its value on $\partial D_0$, and on the latter it is at least that value. Therefore
\begin{equation}\label{eq:Phi}
 \area(G_w(D))\le\Phi(A),\qquad
 \Phi(A):=2\pi\left(1-\frac1{\sqrt{1+A/\pi}}\right),
\end{equation}
where the formula follows by integrating the density over $D_0$ in polar coordinates. The function $\Phi$ is increasing. Since $ah\le(a^2+h^2)/2\le1/2$, the ellipse has planar area at most $\pi/2$, proving~\eqref{eq:ellipse-bound}. The strict comparison there is equivalent to $2/3>16/25$.

The planar hull of five points of the ellipse is a convex polygon with at most five vertices; it may be a segment or a point. A polygon with $k\le5$ vertices in the unit disk has area at most $k/2$: for vertices $z_1,\ldots,z_k$ in boundary order, the shoelace formula gives
\[
 \frac12\left|\sum_{r=1}^k\det(z_r,z_{r+1})\right|\le\frac k2,
 \qquad z_{k+1}=z_1,
\]
where $\det(z_r,z_{r+1})$ is the $2\times2$ determinant with those columns. The linear map taking the unit disk to the ellipse multiplies area by $ah$. Hence the planar hull has area at most
\[
 \frac52ah\le\frac54<\frac{11\pi}{25}.
\]
Gnomonic coordinates preserve convex hulls, so~\eqref{eq:Phi} gives
\[
 \area\bigl(\sconv\{u_1,\ldots,u_5\}\bigr)
 \le\Phi(5/4)<\Phi(11\pi/25)=\frac\pi3. \qedhere
\]
\end{proof}

A member of a finite cover is called \emph{essential} if deleting it destroys the cover. Equivalently, it has a \emph{private point}, a point which belongs to that member and to no other member.

\begin{lemma}\label{lem:essential}
If six sets $E_M$ from Lemma~\ref{lem:sections} cover $\Sph^2$, then each is the union of two spherical ellipses, and every ellipse is essential in the resulting cover.
\end{lemma}

\begin{proof}
If one $E_M$ were a great circle, the remaining at most ten ellipses would have total area at most $10B_0<4\pi$. Adding finitely many great circles, each of area zero, cannot cover the sphere. Hence every $E_M$ is the union of two spherical ellipses.

Suppose one component $C$ were redundant, meaning that the remaining components still covered the sphere. Antipodality implies that $-C$ is also redundant. Since $C\cap(-C)=\varnothing$, every point in either component is covered by a component outside that pair. Thus both may be deleted simultaneously, leaving ten ellipses which cover the sphere. Their total area is less than $4\pi$, a contradiction.
\end{proof}

\section{A covering lemma on the sphere}\label{sec:structure}
A \emph{spherical convex polygon} is the gnomonic image of a closed bounded planar convex polygon with nonempty interior. Its sides are short great-circle arcs, that is, geodesic arcs, since gnomonic projection maps line segments to such arcs. Such a polygon, like a spherical ellipse, is a closed spherically convex region with nonempty interior contained in an open hemisphere.

The proof first counts pair and triple intersections using Euler's formula, and then joins points in the triple intersections to obtain twelve pentagons.

\begin{proposition}\label{prop:structure}
Let $C_1,\ldots,C_{12}$ be spherical ellipses or spherical convex polygons which cover $\Sph^2$ and occur in six antipodal pairs. Suppose each $C_i$ is essential and has spherical diameter less than $\pi/2$. Then there are thirty nonempty pair intersections and twenty nonempty triple intersections. Every intersecting pair belongs to exactly two intersecting triples, and every set belongs to exactly five intersecting triples.

Let
\[
 \mathcal T=\{F\subset\{1,\ldots,12\}:|F|=3,
                   \ \bigcap_{i\in F}C_i\ne\varnothing\}.
\]
Choose any point $a_F\in\bigcap_{i\in F}C_i$ for each $F\in\mathcal T$. Then the twelve spherical convex hulls
\begin{equation}\label{eq:five-hulls}
 P_i=\sconv\{a_F:F\in\mathcal T,\ i\in F\},\qquad 1\le i\le12,
\end{equation}
cover $\Sph^2$, and $P_i\subset C_i$. Each is the hull of five points.
\end{proposition}

\subsection{Euler's formula}
We use Euler's formula in the following form. Subdivide the arrangement of the boundaries of the sets $C_i$ at all intersections, treating coincident boundary portions as single arcs. Add auxiliary arcs if necessary so that each remaining open two-dimensional region is homeomorphic to an open planar ball. For a set which is a union of cells of this subdivision, its Euler characteristic is
\[
 \chi=V_0-E_0+F_0,
\]
where $V_0,E_0,F_0$ count its vertices, open edge pieces and open two-dimensional cells. Subdivision does not change this number. A compact convex region with nonempty interior, a segment or a point has Euler characteristic $1$, and $\chi(\Sph^2)=2$. A connected compact planar region with $h$ holes has Euler characteristic $1-h$, where a hole is a bounded connected component of the complement. These identities follow from the planar Euler formula by including the complementary regions; lower-dimensional edges and isolated boundary contacts are counted in the same subdivision.

The boundaries of spherical ellipses and polygons have such a finite arrangement: apart from coincident portions they have only finitely many intersections. Every intersection or union of the sets is a union of cells in a common subdivision. Inclusion--exclusion for $\chi$ therefore follows by applying ordinary finite inclusion--exclusion separately to vertices, edges and two-dimensional cells. In particular, tangencies and point or segment intersections cause no difficulty.

\subsection{Intersection counts}\label{subsec:intersection-counts}
\begin{proof}[Proof of the intersection assertions]
\textit{Restrictions on intersections.}
If $u,v\in C_i$, then $u\cdot v>0$, by the diameter bound. A set cannot meet its antipode or both members of another antipodal pair. Indeed, if $u\in C_i\cap C_j$ and $v\in C_i\cap(-C_j)$, the diameter bounds give both $u\cdot v>0$ and $u\cdot(-v)>0$. Thus every set has at most five neighbours, where a neighbour is another set with which it intersects.

There is no fourfold intersection. If four sets had a common point $u$, their four antipodal pairs would all miss the great circle $u^\perp\cap\Sph^2$. Only four sets would remain to cover it. By convexity and the diameter bound, each meets this circle in an empty set or an arc of length less than $\pi/2$. These four arcs cannot cover a circle of length $2\pi$.

We also need a restriction on triples of \emph{antipodal-pair labels}. Temporarily label the six pairs by $1,\ldots,6$. Call a triple of pair labels \emph{available} if some choice of one set from each of those pairs has nonempty intersection. Complementary triples of pair labels cannot both be available. Otherwise, if their respective common points were $u$ and $v$, a unit vector perpendicular to both would be missed by all six pairs, contradicting the cover.

\medskip
\noindent\textit{The local and global counts.}
Let $e$ and $t$ be the numbers of nonempty pair and triple intersections of the twelve sets. Every such intersection is convex and has Euler characteristic $1$. Since there are no fourfold intersections, inclusion--exclusion gives
\begin{equation}\label{eq:global-Euler}
 2=\chi(\Sph^2)=12-e+t.
\end{equation}
For a fixed $i$, let $d_i$ be the number of neighbours of $C_i$, and let $t_i$ be the number of triple intersections containing $C_i$. Put
\[
 W_i=C_i\cap\bigcup_{j\ne i}C_j.
\]
The other closed sets cover $\partial C_i$. To see this, approach any boundary point from outside $C_i$ and use the closedness of the finite union of the other sets.

Every nonempty overlap $C_i\cap C_j$ meets $\partial C_i$. Indeed, $C_j$ has a private point outside $C_i$, and the short arc in $C_j$ from an overlap point to this private point reaches $\partial C_i$. Therefore $W_i$ is connected: each convex overlap piece meets the connected set $\partial C_i\subset W_i$.

A private point of $C_i$ is an interior point. In fact, a neighbourhood of it misses all the other closed sets and so must be contained in $C_i$. In gnomonic coordinates, $\partial C_i\subset W_i$ surrounds this point. Thus $W_i$ has at least one hole. Inclusion--exclusion for the overlap pieces $C_i\cap C_j$ yields
\begin{equation}\label{eq:local-Euler}
 d_i-t_i=\chi(W_i)\le0.
\end{equation}
There are no triple intersections among these overlap pieces, since such an intersection would lie in four of the original sets.

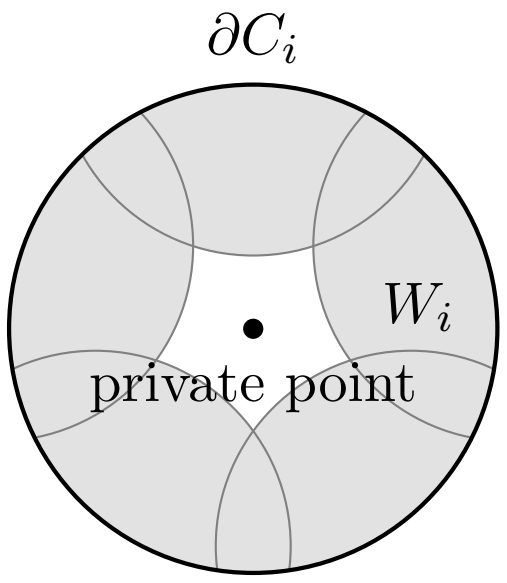
\begin{figure}[ht]
\centering
\begin{tikzpicture}[font=\small,scale=.95]
 \begin{scope}
  \clip (0,0) circle (1.7);
  \foreach \a in {90,18,306,234,162}{\fill[gray!23] (\a:1.87) circle (1.36);}
  \foreach \a in {90,18,306,234,162}{\draw[gray] (\a:1.87) circle (1.36);}
 \end{scope}
 \draw[thick] (0,0) circle (1.7);
 \fill (0,0) circle(2pt);
 \node[below] at (0,-.08) {private point};
 \node at (1.15,.15) {$W_i$};
 \node[above] at (0,1.7) {$\partial C_i$};
 \node[anchor=west,align=left,text width=7.6cm] at (2.4,0)
 {The overlap pieces cover $\partial C_i$ and each reaches that boundary.
 Their union is connected and has a hole.\par\medskip
 Consequently $\chi(W_i)=d_i-t_i\le0$.};
\end{tikzpicture}
\caption{The local Euler count in a gnomonic chart. The arrangement is schematic; only connectedness and the existence of a hole are used.}
\label{fig:local-Euler}
\end{figure}

Since $\sum_i d_i=2e$, $\sum_i t_i=3t$ and $d_i\le5$, equations~\eqref{eq:global-Euler} and~\eqref{eq:local-Euler} imply
\begin{equation}\label{eq:forced-equality}
 2=\sum_{i=1}^{12}\left(1-\frac{d_i}{2}+\frac{t_i}{3}\right)
 \ge\sum_{i=1}^{12}\left(1-\frac{d_i}{6}\right)\ge2.
\end{equation}
Equality throughout forces $d_i=t_i=5$ for every $i$. In particular, $e=30$ and $t=20$.

\medskip
\noindent\textit{Each pair belongs to two triples.}
After forgetting the choice of member within each antipodal pair, the twenty actual triples give ten available triples of pair labels, with each pair label occurring five times. To verify the factor of two, fix one member of one pair. The member of any other pair meeting it is unique, since a set cannot meet both antipodes. Hence each available pair-label triple corresponds to exactly two antipodal actual triples.

The twenty three-element subsets of $\{1,\ldots,6\}$ form ten complementary pairs. Exactly one triple from each complementary pair is available, by the restriction above and the count of ten available triples. Fix two pair labels $a,b$, and let $m$ be the number of available triples containing both. There are four possible triples containing both and four containing neither, paired by complementation. Thus $4-m$ available triples contain neither, and
\[
 10-m-(4-m)=6
\]
available triples contain exactly one of $a,b$. Counting occurrences of $a$ and $b$ gives $5+5=2m+6$, so $m=2$.

It follows that each intersecting pair of actual sets belongs to exactly two actual triples. Indeed, after one member of the first pair is fixed, the member of the second pair in any available triple is forced to be its unique neighbour from that pair.

For each set $C_i$, consider the finite simple graph $L_i$ whose vertices are its five neighbours, and in which distinct neighbours $j,k$ are joined by an edge precisely when $C_i\cap C_j\cap C_k\ne\varnothing$. Here a simple graph has at most one edge between each pair of distinct vertices and has no loops. The degree of a vertex is the number of incident edges. The preceding pair count shows that every vertex of $L_i$ has degree two. A finite simple graph with all degrees two is a disjoint union of cycles, each with at least three vertices. Since $L_i$ has only five vertices, it is a single five-cycle.
\end{proof}

\subsection{The pentagonal regions}\label{subsec:pentagons}
\begin{proof}[Completion of the proof of Proposition~\ref{prop:structure}]
\textit{The embedded graph.}
Choose the points $a_F$ as in the statement. They are distinct, since a common point for two different triples would belong to at least four sets. For each intersecting pair $\{i,j\}$, its two triples have the form $\{i,j,k\}$ and $\{i,j,\ell\}$. Join their chosen points by the short arc
\[
 \gamma_{ij}=[a_{\{i,j,k\}},a_{\{i,j,\ell\}}].
\]
The notation $[u,v]$ here denotes the underlying short arc, without an orientation. Convexity gives $\gamma_{ij}\subset C_i\cap C_j$.

These arcs form an embedded graph $G$ on the sphere: distinct edges meet only at a common endpoint, and no nonincident vertex lies on an edge. Here vertices are the twenty points $a_F$, and edges are the thirty arcs $\gamma_{ij}$. To check the assertion, arcs associated with disjoint pairs of labels cannot meet, because any meeting point would belong to four sets. If $\gamma_{ij}$ meets $\gamma_{ik}$, then $C_i\cap C_j\cap C_k\ne\varnothing$, so $\{i,j,k\}\in\mathcal T$ and both arcs already have $a_{\{i,j,k\}}$ as an endpoint. If their great circles are distinct, their only other possible meeting point is its antipode, which is excluded by containment in $C_i$.

If the arcs lie on the same great circle and overlap beyond their common endpoint, they initially follow the same direction. The nearer of their other endpoints, say $a_{\{i,j,\ell\}}$, would then lie on $\gamma_{ik}$ and hence in four sets. If instead they initially follow opposite directions, they meet only at their common endpoint, since both arcs have length less than $\pi/2$. Finally, a chosen triple point on a nonincident edge would belong to the union of its three labels and the edge's two labels, again at least four sets. Thus $G$ is embedded, also when some of the arcs lie on the same great circle.

\medskip
\noindent\textit{The twelve pentagonal faces.}
The five triple points whose labels contain $i$, together with the five arcs $\gamma_{ij}$, form a simple five-cycle $Q_i$. This follows from the five-cycle $L_i$: its edges correspond to the triple points, and consecutive edges of $L_i$ give consecutive vertices of $Q_i$. A simple cycle is a closed edge path with no repeated vertices other than its first and last vertex.

All of $Q_i$ lies in $C_i$. In gnomonic coordinates, let $D_i$ be the closed bounded polygonal region enclosed by this cycle. It lies in the planar convex hull of its five vertices, so, returning to the sphere,
\begin{equation}\label{eq:regions}
 D_i\subset P_i\subset C_i.
\end{equation}
No graph vertex outside $Q_i$ lies in $C_i$: its defining triple omits $i$, so membership in $C_i$ would give a fourfold intersection. Nor can an edge be a chord of $Q_i$, meaning an edge between two of its vertices which is not a cycle edge. Two distinct triples containing $i$ can share a pair only if that pair contains $i$, and the corresponding edge is already in $Q_i$.

Consequently no graph edge enters the interior of $D_i$. Such an edge would have to cross the cycle, have a vertex inside it, or be a chord, and all three possibilities have been excluded. Thus the interior of $D_i$ is a face of $G$, where a face means a connected component of $\Sph^2\setminus G$.

The intersection graph, whose vertices are the sets and whose edges join intersecting pairs, is connected: otherwise the union of sets in one of its components and the union of all the others would be disjoint nonempty closed sets covering the connected sphere. Cycles $Q_i,Q_j$ for intersecting sets share their edge $\gamma_{ij}$, so $G$, the union of these cycles, is connected as well. Euler's formula therefore gives
\[
 \#\{\text{faces of }G\}=30-20+2=12.
\]
The twelve $D_i$ are distinct. Indeed, two different sets have at most two triple points in common, so they cannot supply the same five-cycle. The $D_i$ therefore account for all faces. Their union, including their boundaries, is $\Sph^2$. Inclusion~\eqref{eq:regions} proves that the hulls $P_i$ cover the sphere.
\end{proof}

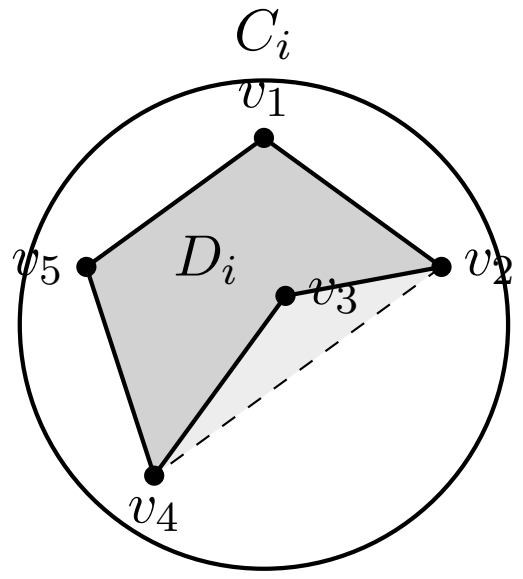
\begin{figure}[ht]
\centering
\begin{tikzpicture}[font=\small,scale=.95]
 \coordinate (v1) at (90:1.3);\coordinate (v2) at (18:1.3);
 \coordinate (v3) at (.15,.20);\coordinate (v4) at (234:1.3);
 \coordinate (v5) at (162:1.3);
 \draw[thick] (0,0) circle (1.7);
 \fill[gray!14] (v1)--(v2)--(v4)--(v5)--cycle;
 \fill[gray!35] (v1)--(v2)--(v3)--(v4)--(v5)--cycle;
 \draw[thick] (v1)--(v2)--(v3)--(v4)--(v5)--cycle;
 \draw[dashed] (v2)--(v4);
 \foreach \r in {1,...,5}{\fill (v\r) circle(2pt);}
 \node[above] at (v1) {$v_1$};\node[right] at (v2) {$v_2$};
 \node[right] at (v3) {$v_3$};\node[below] at (v4) {$v_4$};
 \node[left] at (v5) {$v_5$};\node at (-.4,.45) {$D_i$};
 \node[above] at (0,1.7) {$C_i$};
 \node[anchor=west,align=left,text width=7.6cm] at (2.4,0)
 {The five arcs enclose the face $D_i$ (darker shading).
 The convex hull $P_i$ consists of both shaded regions.\par\medskip
 Both lie in $C_i$.};
\end{tikzpicture}
\caption{A schematic gnomonic view of one pentagonal face. For the cyclically ordered neighbours $j_1,\ldots,j_5$, put $v_r=a_{\{i,j_r,j_{r+1}\}}$, reading indices modulo five. The edge from $v_{r-1}$ to $v_r$ is $\gamma_{ij_r}$. The regions $D_i$ have disjoint interiors, but their convex hulls $P_i$ may overlap.}
\label{fig:pentagon}
\end{figure}

\section{An alternative topological proof}\label{sec:topological}
We give another proof of the intersection assertions in
Proposition~\ref{prop:structure}. The nerve theorem transfers the problem to
a simplicial complex on six vertices. A classification of such complexes
then determines all pair and triple intersections. The subsequent geometric
construction of the covering pentagons is the same as in Section~\ref{subsec:pentagons}.

\begin{proof}
Assume the hypotheses of Proposition~\ref{prop:structure}, and write
$C_i^+,C_i^-=-C_i^+$, $1\le i\le6$, for the antipodal pairs.
Let $q:\Sph^2\to\RP^{2}$ be the quotient map.

\medskip
\noindent\textit{An open cover of the projective plane.}
In a gnomonic chart for each $C_i^+$, add a sufficiently small closed disk
to its planar image and take the interior of the resulting convex set.
Use the antipodal chart for $C_i^-$.
This gives antipodal open spherically convex sets
$U_i^\pm\supset C_i^\pm$ of diameter less than $\pi/2$ with the same
nonempty intersections as the original sets.
Such a choice is possible because the diameter bounds are strict, and if
an empty intersection became nonempty under
arbitrarily small enlargements, compactness would give a point in the
original intersection. A single enlargement suffices because there are
only finitely many subfamilies.

The six sets $\mathcal U_i=q(U_i^+)=q(U_i^-)$ form an open cover of $\RP^{2}$.
Let $J$ be a nonempty set of indices for which
$\bigcap_{i\in J}\mathcal U_i\ne\varnothing$, and fix $i_0\in J$.
For each other $j\in J$, precisely one of $U_j^+$ and $U_j^-$ meets
$U_{i_0}^+$: at least one must, and the diameter bound excludes both,
as in the first paragraph of Section~\ref{subsec:intersection-counts}.
With the corresponding signs $\sigma_j$, the quotient map gives a homeomorphism
\[
 U_{i_0}^+\cap\bigcap_{j\in J\setminus\{i_0\}}U_j^{\sigma_j}
 \xrightarrow{\ q\ }
 \bigcap_{i\in J}\mathcal U_i.
\]
The set on the left is convex in a gnomonic chart for $U_{i_0}^+$,
so it is contractible. Thus the quotient cover is a \emph{good cover}:
every nonempty intersection of its members is contractible.

\medskip
\noindent\textit{The six-vertex nerve.}
Let $\mathcal N$ be the nerve of this cover: its vertices are $1,\ldots,6$,
and a set of vertices forms a simplex when the corresponding $\mathcal U_i$
have nonempty intersection. Write $|\mathcal N|$ for its geometric realization.
The nerve theorem~\cite{Hatcher}*{Corollary~4G.3} gives
\[
 |\mathcal N|\simeq\RP^{2},\qquad
 \pi_1(|\mathcal N|)\cong\mathbb Z/2\mathbb Z.
\]
Govc, Marzantowicz, Michalak and Pave\v{s}i\'{c}~\cite{GMMM}*{Theorem~3.6
and Corollary~3.7} show that a connected simplicial complex on at most six
vertices with this fundamental group is the unique six-vertex triangulation
of $\RP^{2}$. Hence $\mathcal N$ is that triangulation, not just homotopy
equivalent to it.

If $e$ and $f$ count its edges and triangles, then $6-e+f=1$ and $3f=2e$,
so $e=15$ and $f=10$. Every pair of vertices is therefore joined by an edge.
The link of a vertex has its five neighbours as vertices, with an edge for
every triangle containing the original vertex. Since $\mathcal N$
triangulates a closed surface, this link is a cycle, hence a five-cycle.

\medskip
\noindent\textit{Lifting the intersection structure to the sphere.}
Each nonempty simplex of $\mathcal N$ has exactly two lifts to the nerve of
the twelve sets $U_i^\pm$, exchanged by antipodality: fixing one lifted
vertex fixes all signs. In particular, the projection to $\mathcal N$ is
an isomorphism on each vertex link.
Since enlargement preserves intersections, the nerve of the original
twelve sets has twelve vertices, thirty edges, twenty triangles, no
higher-dimensional simplices, and five-cycle links.
Every intersecting pair therefore belongs to two triples, and every
set belongs to five triples. These are the intersection assertions of
Proposition~\ref{prop:structure}.
\end{proof}

To obtain the remaining conclusion of Proposition~\ref{prop:structure},
choose a point in each triple intersection and join the two
points associated with each intersecting pair by their short arc.
The five-cycle links give a pentagonal circuit in each set.
The embedded-graph and face arguments in Section~\ref{subsec:pentagons} show that the
twelve enclosed regions cover $\Sph^2$, each lying in the convex hull of
its five vertices. For the ellipses arising from balls of radius $1/2$,
Lemma~\ref{lem:area} bounds each hull's area by a number strictly less than
$\pi/3$. This gives the same contradiction as in the proof of
Theorem~\ref{thm:main}.

\section{Completion of the proof}\label{sec:completion}
\begin{proof}[Proof of Theorem~\ref{thm:main}]
Compactness, diameter and the convex hull assertions were established in Section~\ref{sec:sections}. Suppose, to the contrary, that six balls of radius $1/2$ cover $V$. By Lemma~\ref{lem:sections}, six sets $E_M$ cover $\Sph^2$. Lemma~\ref{lem:essential} shows that they consist of twelve spherical ellipses, each essential in the cover. They occur in antipodal pairs and have diameter less than $\pi/2$.

Proposition~\ref{prop:structure} gives twelve hulls $P_i$, each spanned by five points of its corresponding ellipse, whose union is the sphere. By Lemma~\ref{lem:area}, every $P_i$ has area less than $\pi/3$. Consequently,
\[
 4\pi=\area(\Sph^2)
 \le\sum_{i=1}^{12}\area(P_i)
 <12\cdot\frac\pi3=4\pi,
\]
the required contradiction.
\end{proof}

To obtain the finite example in Corollary~\ref{cor:finite}, we first show that six balls still cannot cover $V$ when their radius is increased to a suitable $R>0.50001$. Increasing the radius lowers the threshold in the quadratic inequalities defining the spherical inverse images. Near the exceptional eigenvalue configuration of Lemma~\ref{lem:sections}, the resulting sets may spread around a great circle. We therefore separate two cases. In the first, the set lies in a narrow spherical band whose area is too small for it to occur in a six-set cover. In the second, it consists of two antipodal ellipses of diameter less than $\pi/2$. The area estimates for these ellipses and their five-point hulls remain strong enough to apply Proposition~\ref{prop:structure} and obtain the same contradiction.

The positive margin between $R$ and $0.50001$ then permits passage to a finite set. A cover of a sufficiently fine finite net in $V$ by balls of radius $0.50001$ would cover all of $V$ after increasing the radii by the mesh size, while keeping them below $R$. Finally, including two points at distance $1$ ensures that the finite set has diameter exactly $1$.

\begin{proof}[Proof of Corollary~\ref{cor:finite}]
Put
\begin{equation}\label{eq:threshold-gap}
 d=\frac1{40000},\qquad t=\tau-d,\qquad
 R=\sqrt{\frac14+\tau d-\frac{d^2}{2}}.
\end{equation}
We first show that six balls of radius $R$ cannot cover $V$.
The calculation in~\eqref{eq:ball-section} places the inverse image of any
such ball inside
\[
 E_M(t)=\{u\in\Sph^2:u^TMu\ge t\},\qquad
 \tr M=0,\quad\norm M=1.
\]
Indeed, for a nonzero centre $\rho M$, the threshold is at least
$(\rho^2+1/3-R^2)/(\sqrt2\rho)\ge t$; a ball centred at zero misses $V$.

Let $\lambda_1\ge\lambda_2\ge\lambda_3$ be the eigenvalues of $M$, and put
$k=\sqrt{6\tau d-3d^2}$.
The trace and norm conditions give $\tau\le\lambda_1\le2\tau$ and
$\lambda_1+\lambda_2=(\lambda_1+\sqrt{2-3\lambda_1^2})/2$.
The latter expression decreases with $\lambda_1$, so
$\lambda_1+\lambda_2\ge2t$ holds exactly when $\lambda_1\le t+k$.

Suppose first that $\lambda_1+\lambda_2\ge2t$.
Then $\lambda_3\le-2t$ and $\lambda_1\le t+k$.
In an orthonormal eigenbasis, every $u\in E_M(t)$ satisfies
\[
 u_3^2\le\frac{\lambda_1-t}{\lambda_1-\lambda_3}
       \le\frac{k}{3t+k}.
\]
Thus $E_M(t)$ lies in a spherical band, and
\begin{equation}\label{eq:band-area}
 \area(E_M(t))\le4\pi\sqrt{\frac{k}{3t+k}}<\frac{8\pi}{25}.
\end{equation}

Suppose instead that $\lambda_1+\lambda_2<2t$.
Then $\lambda_1>t>\lambda_2$, and $E_M(t)$ consists of two antipodal
spherical ellipses with gnomonic semiaxes
\[
 a^2=\frac{\lambda_1-t}{t-\lambda_2},\qquad
 h^2=\frac{\lambda_1-t}{t-\lambda_3},\qquad h\le a<1.
\]
Each ellipse therefore has spherical diameter less than $\pi/2$.
Writing $s=\lambda_1-\tau$, we have $k-d<s\le\tau$, and the trace and
norm conditions yield
\begin{equation}\label{eq:quantitative-axes}
 (ah)^2=\frac{(s+d)^2}{(s-d)(3\tau+s)+d^2}
 \le\left(\frac{\tau+d}{2\tau-d}\right)^2.
\end{equation}
For the inequality, the quotient increases with $s$ on this interval:
its derivative has the sign of $(\tau-d)s-3\tau d+d^2$, which is positive
because $s>k-d>4d$.
Put $m=(\tau+d)/(2\tau-d)$.
The proof of Lemma~\ref{lem:area} bounds the planar areas of the ellipse
and a five-point hull by $\pi ah$ and $5ah/2$, respectively. Hence
\begin{equation}\label{eq:quantitative-area}
 \area(C)\le\Phi(\pi m)<\frac{46\pi}{125},\qquad
 \area\bigl(\sconv\{v_1,\ldots,v_5\}\bigr)
 \le\Phi(5m/2)<\frac\pi3
\end{equation}
for either ellipse $C$ and any $v_1,\ldots,v_5\in C$.

Suppose now that six sets $E_M(t)$ cover $\Sph^2$.
Each has area less than $92\pi/125$. None can satisfy~\eqref{eq:band-area},
since then their total area would be less than
\[
 \frac{8\pi}{25}+5\frac{92\pi}{125}=4\pi.
\]
Thus all six are unions of two ellipses of diameter less than $\pi/2$.
Every ellipse is essential, as in Lemma~\ref{lem:essential}: deleting a
redundant antipodal pair would leave ten ellipses of total area less than
$10(46\pi/125)<4\pi$.
Proposition~\ref{prop:structure} gives twelve five-point hulls covering
$\Sph^2$, contradicting~\eqref{eq:quantitative-area}.
Thus six balls of radius $R$ cannot cover $V$.

Since $\tau>0.408$, we have
\[
 R^2-(0.50001)^2
 =d\left(\tau-\frac25\right)-\frac{d^2}{2}-10^{-10}>0.
\]
Choose $0<\delta<R-0.50001$ and a finite $\delta$-net $F\subset V$.
Include $p(e_1)$ and $p(e_2)$ for two orthogonal unit vectors, so that
$\diam F=1$. If six balls of radius $0.50001$ covered $F$, the balls
with the same centres and radius $0.50001+\delta<R$ would cover $V$,
a contradiction.
\end{proof}

\section{A Borsuk partition of the example}\label{sec:Borsuk}

The regular dodecahedron gives a partition of $V$ into six sets of diameter less than $1$. Thus this example distinguishes Gr\"unbaum's covering problem from Borsuk's partition problem.

\begin{proposition}\label{prop:Borsuk}
The set $V$ can be partitioned into at most six subsets of diameter at most $2\sqrt2/3$.
\end{proposition}

\begin{proof}
Let $\varphi=(1+\sqrt5)/2$ and take the regular dodecahedron with vertices
\[
 (\pm1,\pm1,\pm1),\quad
 (0,\pm\varphi^{-1},\pm\varphi),\quad
 (\pm\varphi^{-1},\pm\varphi,0),\quad
 (\pm\varphi,0,\pm\varphi^{-1}),
\]
where the displayed signs are independent. Each vertex has norm $\sqrt3$. Its twelve faces, projected radially by $x\mapsto x/\norm x$, cover $\Sph^2$ and occur in six antipodal pairs.

For completeness, the face with outward normal $(0,\varphi,1)$ has vertices
\[
 (1,1,1),\ (-1,1,1),\ (0,\varphi^{-1},\varphi),\
 (\varphi^{-1},\varphi,0),\ (-\varphi^{-1},\varphi,0),
\]
all satisfying $(0,\varphi,1)\cdot x=\varphi^2$. Scalar products of distinct vertices in this face are $1$ or $\sqrt5$. The same holds for every face by symmetry. Thus the unit vertex vectors $v_1,\ldots,v_5$ of any radial face satisfy $v_r\cdot v_s\ge1/3$ for all $r,s$.

Any two points of that radial face have the form $u=x/\norm x$ and $v=y/\norm y$, where $x$ and $y$ are convex combinations of $v_1,\ldots,v_5$. Hence $x\cdot y\ge1/3$ and $\norm x,\norm y\le1$, giving $u\cdot v\ge1/3$. Equation~\eqref{eq:distance} now gives
\[
 \norm{p(u)-p(v)}\le\sqrt{1-\frac19}=\frac{2\sqrt2}{3}.
\]
Antipodal faces have the same image under $p$, so the twelve faces give six covering sets of $V$ with the required diameter bound. Assign a point in several sets to the first one in a fixed ordering. This removes overlaps without increasing any diameter and gives the desired partition.
\end{proof}

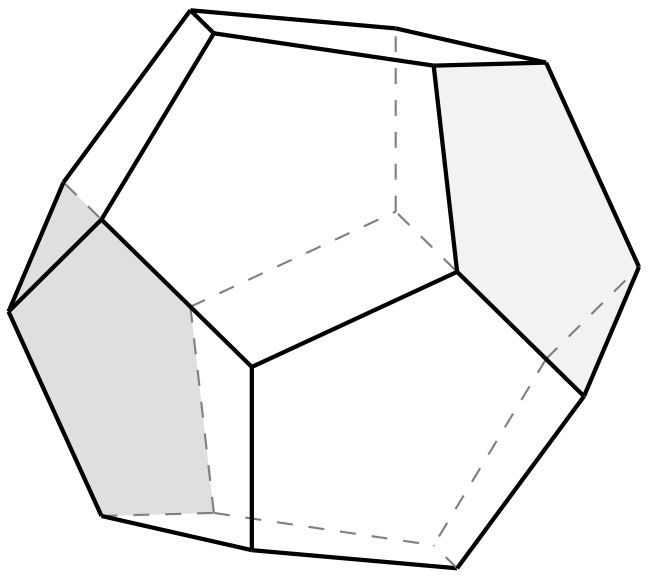
\begin{figure}[ht]
\centering
\begin{tikzpicture}[font=\small]
\coordinate (d0) at (-0.88218,-0.11534);
\coordinate (d1) at (-0.88218,1.84543);
\coordinate (d2) at (1.47029,-0.46136);
\coordinate (d3) at (1.47029,1.49941);
\coordinate (d4) at (-1.47029,-1.49941);
\coordinate (d5) at (-1.47029,0.46136);
\coordinate (d6) at (0.88218,-1.84543);
\coordinate (d7) at (0.88218,0.11534);
\coordinate (d8) at (-0.72695,-1.47937);
\coordinate (d9) at (-1.72145,0.70764);
\coordinate (d10) at (0.47580,0.51383);
\coordinate (d11) at (-0.72695,1.69322);
\coordinate (d12) at (2.08493,0.14777);
\coordinate (d13) at (0.47580,1.72565);
\coordinate (d14) at (0.72695,-1.69322);
\coordinate (d15) at (-2.08493,-0.14777);
\coordinate (d16) at (-0.47580,-1.72565);
\coordinate (d17) at (0.72695,1.47937);
\coordinate (d18) at (1.72145,-0.70764);
\coordinate (d19) at (-0.47580,-0.51383);
\fill[gray!10] (d18)--(d12)--(d3)--(d17)--(d7)--cycle;
\fill[gray!25] (d15)--(d9)--(d0)--(d8)--(d4)--cycle;
\draw[gray,dashed] (d10)--(d13);
\draw[gray,dashed] (d0)--(d10);
\draw[gray,dashed] (d0)--(d9);
\draw[thick] (d1)--(d9);
\draw[thick] (d1)--(d13);
\draw[thick] (d3)--(d12);
\draw[gray,dashed] (d2)--(d12);
\draw[gray,dashed] (d2)--(d10);
\draw[thick] (d3)--(d13);
\draw[thick] (d5)--(d11);
\draw[thick] (d1)--(d11);
\draw[thick] (d9)--(d15);
\draw[thick] (d5)--(d15);
\draw[gray,dashed] (d0)--(d8);
\draw[gray,dashed] (d4)--(d8);
\draw[thick] (d4)--(d15);
\draw[thick] (d6)--(d16);
\draw[thick] (d4)--(d16);
\draw[gray,dashed] (d8)--(d14);
\draw[gray,dashed] (d6)--(d14);
\draw[gray,dashed] (d2)--(d14);
\draw[thick] (d16)--(d19);
\draw[thick] (d5)--(d19);
\draw[thick] (d12)--(d18);
\draw[thick] (d6)--(d18);
\draw[thick] (d7)--(d18);
\draw[thick] (d7)--(d19);
\draw[thick] (d11)--(d17);
\draw[thick] (d3)--(d17);
\draw[thick] (d7)--(d17);
\node[anchor=west,align=left,text width=7.8cm] at (2.9,0)
{Radial projection gives twelve spherical pentagons. Opposite faces have the same image under $p$.\par\medskip
Each pair gives one of six subsets of $V$ with diameter at most $2\sqrt2/3$.};
\end{tikzpicture}
\caption{A regular dodecahedron, with one antipodal pair of faces shaded. Dashed edges lie at the back. The fixed regular faces used here describe a Borsuk partition of $V$; the faces in Proposition~\ref{prop:structure} instead depend on a hypothetical ball cover.}
\label{fig:dodecahedron}
\end{figure}

\section*{AI use disclosure}
The example presented in this paper was found using generative AI tools. The authors verified the proofs and revised the exposition to include both the elementary and the topological approaches.

\begingroup
\linespread{1}\selectfont

\endgroup
\end{document}